\documentclass[12pt]{amsart}
\usepackage{amssymb}
\usepackage{times}
\usepackage{graphicx}

\theoremstyle{plain}

\numberwithin{equation}{section}

\newcommand{\calB}{\mathcal{B}}

\newcommand{\calD}{\mathcal{D}}
\newcommand{\calE}{\mathcal{E}}

\newcommand{\calH}{\mathcal{H}}

\newcommand{\bbZ}{\mathbb{Z}}

\newcommand{\la}{\langle}
\newcommand{\ra}{\rangle}

\def\SL{{\text{SL}}}
\def\Aut{{\text{Aut}}}

\def\Pic{{\text{Pic}}}
\def\Ker{{\text{Ker}}}
\def\O{{\text{O}}}

\def\mod{{\text{mod}}}
\def\dim{{\text{dim}}}
\def\trace{{\text{trace}}}
\def\Hom{{\text{Hom}}}

\begin{document}
\title [Segre cubic] {The Segre cubic and Borcherds products}
\author{Shigeyuki Kond\=o}
\address{Graduate School of Mathematics, Nagoya University, Nagoya,
464-8602, Japan}
\email{kondo@math.nagoya-u.ac.jp}
\thanks{Partially supported by JSPS Grant-in-Aid (S), No 22224001, No 19104001, Japan}

\begin{abstract}
We shall construct a 5-dimensional linear system of holomorphic automorphic forms on a 3-dimensional complex ball by applying
Borcherds theory of automorphic forms.  We shall show that this linear system gives the dual map from the Segre cubic
3-fold to the Igusa quartic 3-fold.
\end{abstract}
\maketitle

\section{Introduction}

The main purpose of this note is to give an application of the theory of automorphic forms on bounded symmetric domains 
of type IV due to Borcherds \cite{B1}, \cite{B2}.  We consider the Segre cubic 3-fold $X$ which is a hypersurface of ${\bf P}^4$
of degree 3 with ten nodes.  The symmetry group ${\mathfrak S}_6$ of degree 6 acts on $X$ as projective transformations.
It is known that the Segre cubic $X$ is isomorphic to the Satake-Baily-Borel compactification of an arithmetic quotient of a 
3-dimensional complex ball
$\calB$ associated to a hermitian form of signature $(1,3)$ defined over the Eisenstein integers (\cite{H1}, \cite{HW}, \cite{H2}).  
The complex ball $\calB$
can be embeded into a bounded symmetric domain $\calD$ of type IV and of dimension 6.  By applying Borcherds' theory of automorphic forms on 
bounded symmetric domains of type IV,
we can construct a 5-dimensional linear system of holomorphic automorphic forms of weight 6.  We shall show that this linear system gives the dual map
from the Segre cubic $X$ to its dual Igusa quartic 3-fold.

B. van Geemen \cite{vG} and B. Hunt \cite{H1} observed that both the Segre cubic 3-fold and the Igusa quartic 3-fold 
are birational to the moduli of ordered 6 points on the projective line.
By taking the triple cover of ${\bf P}^1$ branched along 6 points we get a curve of genus 4 with an automorphism of order 3.
One can consider that the Segre cubic is a compactification of the moduli of such curves.  On the other hand, by taking the double cover
of ${\bf P}^1$ branched along 6 points, we get a hyperelliptic curve of genus 2.  The Igusa quartic is the Satake compactification of
an arithmetic quotient of the Siegel space of degree 2 (J. Igusa \cite{I}, page 397).  
It is classically known that  the dual of the Segre cubic is isomorphic to the Igusa quartic (Baker \cite{Ba}, Chap. V).
We will give an interpretation of this ball quotient as the moduli space of some $K3$ surfaces with an automorphism of order 3 which are
obtained from 6 points on the projective line.

We use an idea of Allcock, Freitag \cite{AF} to construct a linear system of automorphic forms.  In \cite{AF}, 
they consider a 4-dimensional complex ball defined over the Eisenstein integers and construct a 10-dimensional linear system
of automorphic forms.   An arithmetic quoteint of the 4-dimensional complex ball is birational to the moduli space of marked cubic surfaces.  Our complex ball $\calB$ appears as a subcomplex ball of Allcock, Freitag's one, and hence one can restrict Allcock, Freitag's
linear system to $\calB$.  However in this note, instead of using their linear system, we apply Borcherds theory directly to
our situation and get a linear system on $\calB$.  


The plan of this note is as follows.  In \S 2, we recall the Segre cubic $X$ and some divisors on $X$.
In section 3, we mention the complex ball $\calB$, the bounded symmetric domain $\calD$ and Heegner divisors on them.
In section 4, we recall the Weil representation and calculate its character.  In section 5, we shall show that 
there exist holomorphic automorphic forms on the complex ball $\calB$ of
weight $45$, $5$ with known zeros.  These forms will be used to determine the zeros of a member of a 5-dimensional 
linear system of automorphic forms on $\calB$.  In section 6, we construct a 5-dimensional linear system of automorphic forms
and show that this linear system gives the dual map of the Segre cubic.

\section{The Segre cubic threefold}\label{}
In this note we consider the variety $X$ called the {\it Segre cubic} which is defined by 
$$X : \sum_{i=1}^6 x_i =0, \ \sum_{i=1}^6 x_i^3 =0$$
in ${\bf P}^5$.
Obviously the symmetric group ${\mathfrak S}_6$ of degree 6 acts on $X$ projectively.
$X$ has ten nodes which are ${\mathfrak S}_6$-orbits of
$(1:1:1:-1:-1:-1)$.  A linear section of $X$ given by $x_i+x_j=0$ is the union of three projective planes
given by 
$$x_i+x_j =0, \ x_k+x_l =0, \ x_m+x_n=0;$$
$$x_i+x_j =0, \ x_k+k_m =0, \ x_l+x_n=0;$$
$$x_i+x_j =0, \ x_k+x_n =0, \ x_l+x_m=0,$$
respectively where $\{i,j,k,l,m,n\} = \{ 1,2,3,4,5,6\}$.  
On the other hand a linear section of $X$ given by $x_i-x_j=0$ is an irreducible
cubic surface containing four nodes of $X$.  This irreducible cubic surface with four node is projectively unique and is called {\it Cayley cubic surface}.
Thus we have 15 Cayley cubics and 15 planes on the Segre cubic $X$.  
It is known that the dual of $X$ is a quartic 3-fold $Y$ in ${\bf P}^4$ called the {\it Igusa quartic} (\cite{I}).
The dual map $d : X \to Y$ is defined on $X$ except ten nodes and is birational.  It is 
given by a linear system of quadrics through ten nodes.  For more details of these facts, we refer the
reader to \cite{H2}, Chap. 3.

\section{A complex ball quotient}\label{}

It is known that the Segre cubic $X$ is isomorphic to the Satake-Baily-Borel compactification of 
an arithmetic quotient of a 3-dimensional complex ball by a certain 
arithmetic subgroup (\cite{HW}, Theorem 1; \cite{H2}, Chap. 3, 3.2.3).  In this section we recall this fact.

\subsection{A complex ball}\label{}

Let 
$${\calE} = {\bf Z}[\omega], \ \omega = {-1 + \sqrt{-3} \over 2}$$ 
be the ring of Eisenstein integers.
Consider the hermitian lattice 
$$\Lambda = {\calE}^{1,3} = {\calE} \oplus {\calE} \oplus {\calE}
\oplus {\calE}$$ 
with the hermitian form 
$$h(x, y) =
x_{0} {\bar y}_{0} - x_{1} {\bar y}_{1} - x_{2} {\bar y}_{2} - 
x_{3} {\bar y}_{3}.$$
We denote by $\calB$ the complex ball of dimension 3 defined by
$$\calB = \{ x \in {\bf P}(\Lambda \otimes_{\calE} {\bf C}) \ : \ h(x,x) > 0 \}.$$
Let $\Gamma = \Aut(\Lambda)$ be an arithmetic subgroup of the unitary group $U(3,1; {\bf Q}(\sqrt{-3}))$ with respect to the
hermitian form $h(\ , \ )$.  Obviously $\Gamma$ naturally acts on $\calB$.
Under the isomorphism 
$${\calE}/ \sqrt{-3} {\calE} \simeq {\bf F}_{3},$$
the hermitian form $h$ induces a quadratic form $q$ on $\Lambda / \sqrt{-3}
\Lambda$ over ${\bf F}_{3}$.  We define a subgroup $\Gamma(\sqrt{-3})$ of $\Gamma$ by
$$\Gamma(\sqrt{-3}) = \Ker \{\Gamma \longrightarrow \O(q)\}.$$

Let $L$ be the real lattice corresponding to $\Lambda$ with the symmetric bilinear form 
$$\la x, y\ra = h(x,y) + h(y,x).$$
Then $L \cong A_{2} \oplus 
A_{2}(-1)^{3}$ where $A_{2}$ is a root lattice of rank 2, that is, a positive definite lattice of rank 2 defined by
the matrix 
$\begin{pmatrix}2&-1
\\-1&2
\end{pmatrix}$
and $A_2(-1)$ is the negative definite lattice of rank 2 defined by 
$\begin{pmatrix}-2&1
\\1&-2
\end{pmatrix}.$
The action of $\omega$ on $\Lambda$ induces an isometry $\iota$ of $L$ of order 3 without non-zero fixed points.
We denote by $L^*$ the dual of $L$: $L^*= \Hom (L, {\bf Z})$.
Note that $A_2^*/A_2 \cong {\bf F}_3$.
Let $A_{L} = L^{*}/L \cong ({\bf F}_3)^4$ and let $q_{L} : A_L \to {\bf Q}/2{\bf Z}$ be the discriminant quadratic form of $L$ defined by
$q_L(x + L) = \la x,x \ra + 2{\bf Z}$.  The form $q_L$ coincide with $q$, up to scales, under the isomorphism 
$$\Lambda/\sqrt{-3}\Lambda \cong L^*/L.$$ 
We denote by $\O(L)$ the orthogonal group of $L$ and by $\tilde{\O}(L)$ the kernel of the 
natural map $\O(L) \to \O(q_L)$.  Then
the group $\Gamma$ is naturally isomorphic to the subgroup $\O(L,\iota)$ of $\O(L)$ consisting of isometries commuting with $\iota$.
Under this isomorphism the subgroup $\Gamma(\sqrt{-3})$ corresponds to $\O(L,\iota) \cap \tilde{\O}(L)$.

Conversely, first consider the lattice $L$ with an automorphism $\iota$ of order 3 without non-zero fixed points.
Then we can consider $L$ as ${\bf Z}[\omega]$-module by the action $\omega \cdot x = \iota(x)$.
The hermitian form $h$ is given by
$$h(x,y) = {1\over 2}\{ {\sqrt{-3} \over 3}\la 2\iota(x)+x, y\ra + \la x, y \ra\}.$$
Define
$$\calD =\{ v \in {\bf P}(L\otimes {\bf C}) \ : \ \la v, v \ra =0, \ \la v, \bar{v} \ra > 0 \}.$$
Then $\calD$ is a disjoint union of two copies of a bounded symmetric domain of type IV and of dimension 6.
Consider the action of $\iota$ on $L\otimes {\bf C}$.  Since $\iota$ is defined over ${\bf Z}$ and has no non-zero fixed vectors, 
the eigen-spaces $V_{\omega}, V_{\bar{\omega}}$ with the eigen-value $\omega, \bar{\omega}$ respectively are isomorphic to
${\bf C}^4$.  Moreover the restriction $\la v, \bar{v}\ra$ on $V_{\omega}$ is a hermitian form of signature $(1,3)$.
Note that $\la v, v \ra =0$ for any vector $v$ in $V_{\omega}$ or $V_{\bar{\omega}}$ because $\la v, v\ra = \la \iota(v), \iota(v)\ra =
\la \omega v, \omega v \ra$ or $\la v, v \ra = \la \bar{\omega} v , \bar{\omega} v \ra$.
Let $i : \Lambda \to \Lambda \otimes_{\bf Z} {\bf C} = L\otimes {\bf C}$ be the inclusion map and let $p : L\otimes {\bf C} \to V_{\omega}$
be the projection.  
For any $\xi \in L$, write $\xi = \xi_{\omega} + \xi_{\bar{\omega}}$ as an element in $V_{\omega} \oplus V_{\bar{\omega}}$.
Then we can easily see that 
$$h(\xi, \xi) = \la \xi_{\omega}, \bar{\xi_{\omega}} \ra.$$
Hence the map $p\circ i : \Lambda \to V_{\omega}$ is an isometry  which induces an isomorphism from
$\calB$ to the subdomain 
$$\calD \cap {\bf P}(V_{\omega}) = \{ v \in {\bf P}(V_{\omega}) \ : \ \la v, \bar{v} \ra > 0 \}$$ 
of $\calD$.  Thus the complex ball $\calB$ can be embedded into $\calD$. 

\subsection{Roots and reflections}\label{}
Following to \cite{AF}, we recall roots and reflections of the hermitian lattice $\Lambda$.
A vector $a \in \Lambda$ is called a {\it short root} (resp. {\it long root}) if $h(a,a)=-1$ (resp. $h(a,a)=-2$).  For a short root or long root
$a$, consider the following isometry $r_{a,\zeta}$ of $\Lambda$ with respect to $h$:
$$r_{a,\zeta} : x \to x - (1-\zeta){h(a,v) \over h(a,a)}a.$$
If $a$ is a short root and $\zeta$ is a primitive third root of unity $\omega$, $r_{a,\omega}$ is an isometry of $\Lambda$ of order three  sending $a$ to $\omega a$.  We call $r_{a,\omega}$ a {\it trireflection}.   If $a$ is a short root or long root, and $\zeta =-1$, 
then $r_{a,-1}$ is a {\it reflection} in $\Gamma$ which is an isometry of order two sending $a$ to $-a$.  

For a short root $a$ in $\Lambda$, denote by $r$ the corresponding $(-2)$-vector in $L$.  Then
the trireflection $r_{a,\omega}$ induces an isometry 
$$s_r \circ s_{\iota(r)} : x \to x + \la x, r \ra r + \la x, \iota(r)\ra r + \la x, \iota(r)\ra \iota(r)$$
of $L$ where $s_r : x \to x+\la x, r\ra r$ is a reflection associated to $r$.  
On the other hand, $r_{a,-1}$ induces an isometry of $L$:
$$x \to x + 2\la {r + 2\iota(r) \over 3}, x\ra \iota(r) + 2\la {2r +\iota(r)\over 3}, x\ra r.$$

For a long root $a \in \Lambda$, denote by $r$ the corresponding $(-4)$-vector in $L$.
Then $r_{a,-1}$ induces an isometry of $L$:
$$x \to x + \la {r + 2\iota(r) \over 3}, x\ra \iota(r) + \la {2r +\iota(r)\over 3}, x\ra r.$$

For $a \in \Lambda$, we denote by $\bar a$ the image of $a$ in $\Lambda/\sqrt{-3}\Lambda$.
We call the images of short roots (resp. long roots) in $\Lambda/\sqrt{-3}\Lambda$ the short roots (resp. long roots), too.
We also denote by ${\bar r}_{a,\zeta}$ the isometry on $\Lambda/\sqrt{-3}\Lambda$ induced by $r_{a,\zeta}$.
Note that if $a$ is a short root, then $r_{a,\omega}$ is contained in $\Gamma (\sqrt{-3})$, that is,
${\bar r}_{a,\omega}$ acts trivially on $\Lambda/\sqrt{-3}\Lambda$.  On the other hand, ${\bar r}_{a,-1}$ acts
on $\Lambda/\sqrt{-3}\Lambda$ as a reflection associated to $\bar a$.

\subsection{Lemma}\label{orbits}

(1) {\it The group $\Gamma$ acts transitively on the primitive isotropic vectors, on the short roots and on the long roots, respectively.  }
\smallskip

(2) {\it Let $a_1$, $a_2$ be two isotropic vectors, or two short roots, or two long roots.  Then 
$a_1$ and  $a_2$  are equivalent under $\Gamma(\sqrt{-3})$ if and only if their images 
 in $\Lambda/\sqrt{-3}\Lambda$ coincide.}
\smallskip

(3) {\it The number of non-zero isotropic vectors, short roots or long roots in $\Lambda/\sqrt{-3}\Lambda$ is $20$, $30$ or $30$ respectively.}

\smallskip
(4) {\it The map $\Gamma \to \O(q)$ is surjective and $\Gamma/\Gamma(\sqrt{-3}) \simeq \O(q) \simeq {\mathfrak S}_{6} \times \bbZ/2\bbZ$.  }
\begin{proof}
In case of the hermitian lattice $\calE^{1,4}$, Allcock, Carlson, Toledo proved the same assertion (1) (\cite{ACT}, Theorems 7.21, 11.13),
and Allcock, Freitag (\cite{AF}, Proposition 2.1) proved the assertions (2), (3).  
The same proof works in our case $\calE^{1,3}$.  The last assertion is a well know.  For example,
see \cite{C}, page 4.
\end{proof}

\subsection{Ball quotient and Heegner divisors}\label{}
We denote a vector $\alpha\in A_L = ({\bf F}_3)^4$ by $\alpha=(x_1,x_2,x_3,x_4)$ where $x_i \in {\bf F}_3$ is in the $i$-th factor of 
$L^*/L = A_{2}^*/A_2 \oplus (A_{2}(-1)^*/A_2(-1))^{\oplus 3}.$
Then an elementary calculation shows the following:

\subsection{Lemma}\label{types}
(i) {\it $A_{L}$ consists of the following $81$ vectors}:
\smallskip

Type $(00)$: $0$;
\smallskip

Type $(0)$: $ \alpha \not= 0, q_{L}(\alpha) = 0, \# \alpha = 20,$
\smallskip

$\alpha = (\pm 1, \pm 1, 0, 0), (\pm 1, 0, \pm 1, 0), (\pm 1, 0, 0, \pm 1),
(0, \pm 1, \pm 1, \pm 1);$
\smallskip

Type $(1)$: $q_{L}(\alpha) = -4/3, \# \alpha = 30,$
\smallskip

$\alpha = (\pm 1, 0,0,0), (\pm 1, \pm 1, \pm 1, \pm 1), 
(0, \pm 1, \pm 1, 0), (0, \pm 1, 0, \pm 1), (0, 0, \pm 1, \pm 1);$
\smallskip

Type $(2)$: $q_{L}(\alpha) = -2/3, \# \alpha = 30,$ 
\smallskip

$\alpha = (0, \pm 1, 0, 0), (0, 0, \pm 1, 0), (0, 0, 0, \pm 1),
(\pm 1, \pm 1, \pm 1, 0), (\pm 1, \pm 1, 0, \pm 1), 
(\pm 1, 0, \pm 1, \pm 1).$

\medskip

\subsection{Lemma}\label{}
{\it 
Under the canonical isomorphism $\Lambda/\sqrt{-3}\Lambda \cong A_L=L^*/L$,
the set of short roots $($resp. long roots$)$ in $\Lambda/\sqrt{-3}\Lambda$ correspond to the set of vectors of
norm $-2/3$ $($resp. vectors of norm $-4/3$$)$ in $A_L$.  Also the set of isotropic vectors in $\Lambda/\sqrt{-3}\Lambda$
correspond to the set of isotropic vectors in $A_L$.}

\medskip

Let $\alpha$ be a non-isotropic vector in $A_L$.  
For a given $n \in {\bf Z}$, $n <0$, we consider a Heegner divisor $\calD_{\alpha, n}$ which is the union of 
the orthogonal complements $r^{\perp}$ in $\calD$ where $r$ varies over the vectors in $L^*$ satisfying $\la r, r\ra = n$ and
$r \ \mod \ L = \alpha$.  Obviously $r^{\perp}$ is a bounded symmetric domain of type IV and 
of dimension 5.  In case $q_L(\alpha) = -2/3$ (resp. $q_L(\alpha) = -4/3$) and $n=-2/3$ (resp. $n=-4/3)$,
we denote $\calD_{\alpha, -2/3}$ (resp. $\calD_{\alpha, -4/3}$) by $\calD_{\alpha}$ for simplicity and call it 
{\it $(-2/3)$-Heegner divisor} (resp. $(-4/3)$-{\it Heegner divisor}).

\subsection{Proposition}\label{}
{\it 
The Segre cubic $X$ is isomorphic to the Satake-Baily-Borel compactification $\bar{\calB}/\Gamma(\sqrt{-3})$
of the quotient ${\calB}/\Gamma(\sqrt{-3})$
which is, set theoretically, the union of ${\calB}/\Gamma(\sqrt{-3})$ and ten cusps
corresponding to 
ten non-zero isotoropic vectors in $A_L/\{\pm 1\}$.  These ten cusps correspond to ten nodes of
the Segre cubic $X$.  
}
\begin{proof}
The assertion follows from \cite{H2}, \S 3.2.
\end{proof}
Also $\bar{\calB}/\Gamma(\sqrt{-3})$ contains some divisors called {\it Heegner divisors}.
Let $\alpha$ be a short root in $\Lambda/\sqrt{-3}\Lambda$.  
Let $a$ be a short root in $\Lambda$ with $a \ \mod\ \sqrt{-3}\Lambda = \alpha$.
We denote by $a^{\perp}$ the orthogonal complement of $a$ in $\calB$ which is a complex ball of dimension $2$.
Let 
$$\calH_{\alpha} = \bigcup_{a} \ a^{\perp}$$ 
where $a$ moves on the set of all short roots satisfying $a \ \mod\ \sqrt{-3}\Lambda = \alpha$.
The image of $\calH_{\alpha}$ in $\bar{\calB}/\Gamma(\sqrt{-3})$ is denoted by $\bar{\calH}_{\alpha}$ and is called a
$(-1)$-{\it Heegner divisor}.  There exist 15 $(-1)$-Heegner divisors $\bar{\calH}_{\alpha}$ corresponding to 15 short roots
$\alpha \in (\Lambda/\sqrt{-3}\Lambda)/\{\pm 1\}$, $q(\alpha)=-1$. 

Similarly we can define 15 $(-2)$-Heegner divisors $\bar{\calH}_{\alpha}$ corresponding to 15 long roots 
$\alpha \in (\Lambda/\sqrt{-3}\Lambda)/\{\pm 1\}$, $q(\alpha)=-2$.

Finally we compare the Heegner divisors in $\calD$ and in $\calB$.
Let $r \in L^*$ be a $(-2/3)$- or $(-4/3)$-vector.  Then both $\iota(r)$ and $\iota^2(r)$ are
$(-2/3)$- or $(-4/3)$-vectors and $r \ \mod \ L = \iota(r)\ \mod\ L = \iota^2(r) \ \mod \ L$.  Note that  
$r^{\perp}$, $\iota(r)^{\perp}$ and $\iota^2(r)^{\perp}$ in $\calD$  are different, but their restrictions to $\calB$ are the same.  Thus
we have 

\subsection{Lemma}\label{restriction}
$$\calD_{\alpha} \cap \calB = 3\calH_{\alpha}$$
{\it where we identify $(-2/3)$- (resp. $(-4/3)$-) vectors in
$A_L$ and short roots $($resp. long roots$)$ in $\Lambda/\sqrt{-3}\Lambda$.}

\subsection{Interpretation via $K3$ surfaces}\label{k3}

The complex ball quotient ${\calB}/\Gamma(\sqrt{-3})$ can be considered as the moduli space of lattice polarized $K3$ surfaces.  
The following is essentially given in \cite{DGK}.  Let $N = U\oplus E_6(-1) \oplus A_2(-1)^{\oplus 3}$.  Then $N$ can be primitively embedded
into the $K3$ lattice $M=U^{\oplus 3}\oplus E_8(-1)^{\oplus 2}$ whose orthogonal complement is isomorphic to 
$L=A_2\oplus A_2(-1)^{\oplus 3}$.  Here $A_m, \ E_k$ are positive definite root lattices defined by the Cartan matrix of type $A_m,\ E_k$, and for a lattice
$(L, \la, \ra)$ we denote by $L(-1)$ is the lattice $(L, -\la, \ra)$. 
In the following we consider $N$ and $L$ as sublattices of $M$.
The isometry $\iota$ of $L$ of order 3 acts trivially on $L^*/L$ and hence it can be extended to
an isometry $\tilde{\iota}$ of $M$ acting trivially on $N$.  Let $\omega \in \calB$ with the property $\omega^{\perp} \cap M = N$.
Let $S$ be a $K3$ surface and let $\alpha_S : H^2(S,{\bf Z}) \to M$ be an isometry satisfying 
$(\alpha_S\otimes {\bf C})(\omega_S) = \omega$ where $\omega_S$ is a holomorphic 2-form on $S$.
By definition the Picard lattice of $S$ is isomorphic to $N$.  Note that the isometry $\tilde{\iota}$ preserves $\omega_S$ and acts trivially
on the Picard lattice.  It now follows from the Torelli type theorem for $K3$ surface that $\tilde{\iota}$ can be represented 
by an automorphism $\sigma$ on $S$ of order 3.  Thus an open set of ${\calB}/\Gamma(\sqrt{-3})$  is the moduli of such pairs 
$(S, \sigma)$ of $K3$ surfaces $S$ with an automorphism $\sigma$ of order 3.  

In the following, we shall show that $S$ is canonically obtained from six points on ${\bf P}^1$.
Let $Q={\bf P}^1 \times {\bf P}^1$.  Let $(u_0:u_1), (v_0:v_1)$ be homogeneous coordinates of the first and the second factor of $Q$.
Let $p_1,..., p_6$ be distinct six points on ${\bf P}^1$.  Consider the divisors on $Q$ defined by
$$L_i = {\bf P}^1 \times \{p_i\} \ (1\leq i \leq 6),\ D_0 = \{0\} \times {\bf P}^1, \  D_1=\{1\} \times {\bf P}^1, \ D_{\infty} = \{\infty\} \times {\bf P}^1.$$
Let $\tilde{Q} \to Q$ be the blow ups of the 18 points on $Q$ which are the intersection of $L_1,..., L_6$ and $D_0, D_1, D_{\infty}$.  
We denote by the $\tilde{L}_1,..., \tilde{L}_6, \tilde{D}_0, \tilde{D}_1$ or  $\tilde{D}_{\infty}$ the strict transform of $L_1$,..., $L_6$,
$D_0$, $D_1$ or $D_{\infty}$ respectively.  Let $\pi : \tilde{X} \to \tilde{Q}$ be the triple covering of $\tilde{Q}$ branched along 
$\tilde{L}_1+ \cdot\cdot\cdot + \tilde{L}_6 + \tilde{D}_0 + \tilde{D}_1 + \tilde{D}_{\infty}$.
Then $\pi^{-1}(\tilde{L}_i)$ is a $(-1)$-curve.  Let $\tilde{X} \to S$ be the contraction of $\pi^{-1}(\tilde{L}_i)$ to the points $q_i$.   
We can easily see that $S$ is a K3 surface.
The projection from $Q$ to the second factor ${\bf P}^1$ induces an elliptic fibration $p : S \to {\bf P}^1$ which has six singular
fibers of type IV in the notation of Kodaira and three sections.  Here three components of the singular fiber of type IV over $p_i$
correspond to three exceptional curves over the three intersection points of $L_i$ and $D_0, D_1, D_{\infty}$ and
three sections correspond to $D_0, D_1, D_{\infty}$.  The classes of components of fibers and a section generate a sublattice of
the Picard lattice $\Pic(S)$ isomorphic to $U\oplus A_2(-1)^{\oplus 6}$.  By adding other two sections, we have a sublattice in $\Pic(S)$
isomorphic to $N = U \oplus E_6(-1)\oplus A_2(-1)^{\oplus 3}$.  The covering transformation of $\tilde{S}\to \tilde{Q}$ induces an automorphism
$\sigma$ of $S$ of order 3. Note that the set of fixed points of $\sigma$ consists of six isolated points $q_1,..., q_6$ and three sections.
Since $\sigma$ has a fixed curve as its fixed points, $\sigma^*(\omega_S) = \zeta_3\omega_S$ where $\omega_S$ is a non-zero holomorphic 
2-form on $S$ and $\zeta_3$ is a primitive cube root of unity.  Thus we have a pair $(S, \sigma)$ of a $K3$ surface and an automorphism of 
order 3.  This $K3$ surface appears as a degeneration of $K3$ surfaces associated to a smooth cubic surface given in \cite{DGK}.

Next we consider the case two points among 6 points coincide.  In this case, similarly, we have an elliptic $K3$ surface $S'$ with one singular fiber of
type ${\rm VI}^*$, four singular fibers of type VI and three sections.  The Picard lattice of $S'$ is isomorphic to $U\oplus E_6(-1)^2\oplus A_2(-1)$ and
its transcendebtal lattice is isomorphic to $A_2\oplus A_2(-1)^{\oplus 2}$.  Thus the period domain of $K3$ surfaces $S'$ is a subdomain of $\calB$
the orthogonal complement of $A_2(-1)$, that is, a $(-1)$-Heegner divisor.  Thus we have

\subsection{Proposition}\label{15planes}
{\it $15$ $(-1)$-Heegner divisors bijectively correspond to $15$ planes on the Segre cubic $X$.}
\begin{proof}  
It is known that 15 planes on the Segre cubic correspond to the moduli of 6 points on the projective line in which two points coincide
(\cite{H2}, Proposition 3.2.7).  Hence we have the assertion.
\end{proof}

We shall show that 15 $(-2)$-Heegner divisors correspond to 15 Cayley cubics on the Segre cubic $X$
(see Lemma \ref{Cayley}).  

\section{Weil representation}\label{weil}


In this section, we recall a representation of $\SL(2,{\bf Z})$ on the group ring ${\bf C}[A_L]$ called the Weil representation.
In following Table 1, for each vector $u \in A_L$ of given type, 
$m_j$ is the number of vectors $v$
of the same type with $\langle u, v \rangle = 2j/3$.

\begin{table}[h]
\[
\begin{array}{rlllllllllllllllllllllll}
u& 00&00&00&00&0&0&0&0& 1&1&1&1&2&2&2&2\\
v&00&0&1&2&00&0&1&2&00&0&1&2&00&0&1&2\\
m_0&1&20&30&30&1&2&12&12&1&8&12&6&1&8&6&12\\
m_1&0&0&0&0&0&9&9&9&0&6&9&12&0&6&12&9\\
m_2&0&0&0&0&0&9&9&9&0&6&9&12&0&6&12&9\\
\end{array}
\]
\caption{}
\end{table}

Let
$T =
\begin{pmatrix}1&1
\\0&1
\end{pmatrix},\ 
S =
\begin{pmatrix}0&-1
\\1&0
\end{pmatrix}$
be a generator of $\SL(2,{\bf Z})$.
Let $\rho$ be the Weil representation of $\SL(2, {\bf Z})$ on ${\bf C}[A_L]$ defined by:
$$\rho(T)(e_{\alpha}) = exp(\langle \alpha, \alpha
\rangle /2 ) e_{\alpha}, \
\rho(S)(e_{\alpha}) = {-1 \over \sqrt{\mid A_L\mid}} \sum_{\delta} 
exp(- \langle \delta, \alpha \rangle ) e_{\delta}.$$
Note that  the action $\rho$ factorizes the action of $\SL(2,{\bf Z}/3{\bf Z})$ which is denoted by the same symbol $\rho$.
The conjugate classes of $\SL(2, {\bf Z}/3{\bf Z})$ consist of
$\pm E, S, \pm ST, \pm ST^{2}$. 
Let $\chi_{i}$ $(1 \leq i \leq 7)$ be the characters of irreducible
representations of $\SL(2, {\bf Z}/3{\bf Z})$.  The following Table 2 is the character table of
$\SL(2, {\bf Z}/3{\bf Z})$.  Here $\omega = {-1 + \sqrt{-3} \over 2}$  and the last line means the number of elements in a given conjugate class.

\begin{table}[h]
\[
\begin{array}{rlllllllllllllllllllllll}
 & E&-E&S&ST^2&-ST^2&ST&-ST\\
\chi_1&1&1&1&1&1&1&1\\
\chi_2&3&3&-1&0&0&0&0\\
\chi_3&1&1&1&\omega^2&\omega^2&\omega&\omega\\
\chi_4&1&1&1&\omega&\omega&\omega^2&\omega^2\\
\chi_5&2&-2&0&-\omega&\omega&\omega^2&-\omega^2\\
\chi_6&2&-2&0&-1&1&1&-1\\
\chi_7&2&-2&0&-\omega^2&\omega^2&\omega&-\omega\\
&1&1&6&4&4&4&4\\
\end{array}
\]
\caption{}
\end{table}

\subsection{Lemma}\label{char}
{\it Let $\chi$ be the character of the representation of
$\SL(2, {\bf Z}/3{\bf Z})$ on ${\bf C}[A_{L}]$.
Let $\chi = \sum_{i} m_{i} \chi_{i}$ be
the decomposition into irreducible characters.  Then
$m_{1} = 1, m_{2} = 10, m_{3} = m_{4} = 5, m_{5} = 5, m_{6} = 10,
m_{7} = 5.$}

\begin{proof}
By definition of $\rho$ and the Table 1, we can easily see that $\trace (E) = 3^4$, $\trace (-E) = 1$, $\trace (S) = 1$, $\trace(ST^2)=-9$,
$\trace(-ST^2)= 1$, $\trace(ST)=1$, $\trace(-ST)=-9$.  The assertion now follows from the Table 2.
\end{proof}

\subsection{Definition}\label{5-dim}
Let $V$ be the 5-dimensional subspace of ${\bf C}[A_{L}]$ which is the direct sum of
irreducible representations of $SL(2, {\bf F}_{3})$ with the character $\chi_{3}$ in Lemma \ref{char}.  In \S \ref{additive}, we associate
a 5-dimensional space of automorphic forms on $\calB$ to $V$.  We remark that there is an another 5-dimensional subspace
in ${\bf C}[A_L]$ which is a direct sum of irreducible representations of $SL(2, {\bf F}_{3})$ with the character $\chi_{4}$.  The author
does not know whether this subspace corresponds to an interesting linear system of automorphic forms on $\calB$.

\section{Borcherds products}\label{}

Borcherds products are automorphic forms on $\calD$ 
whose zeros and poles lie on Heegner divisors.
In this section, we shall show that there exist automorphic forms $\Phi_{45}, \Phi_{5}$ on the complex ball $\calB (\subset \calD)$ of
weight $45$, $5$ whose zero divisors are $(-2)$-Heegnear divisor, $(-1)$-
Heegnear divisor respectively.  

To show the existence of such 
Borcherds products, we introduce the {\it obstruction space}
consisting of all vector valued elliptic modular forms $\{ f_{\alpha} \}_{\alpha \in A_{L}}$ of weight $(2+6)/2=4$ and with respect to
the dual representation $\rho^*$ of $\rho$:
$$f_{\alpha}(\tau + 1) = e^{-\pi \sqrt{-1}\ \la \alpha, \alpha \ra}f_{\alpha}(\tau), \quad
f_{\alpha}(-1/\tau) = -{\tau^{4} \over 9} 
\sum_{\beta} e^{2\pi \sqrt{-1} \ \la \alpha, \beta \ra} f_{\beta}(\tau).$$
We shall apply the next theorem to show the existence of 
some Borcherds products.

\subsection{Theorem}\label{freitag}(Borcherds \cite{B2}, Freitag \cite{F}, Theorem 5.2)
{\it A linear combination
$$\sum_{\alpha \in A_{L}, n<0} c_{\alpha, n} \calD_{\alpha, n}, \ c_{\alpha, n} \in {\bf Z}$$
of Heegner divisors
is the divisor of an automorphic form on $\calD$ 
of weight $k$ if for every cusp form 
$$f = \{f_{\alpha}(\tau)\}_{\alpha\in A_{L}}, 
\ f_{\alpha}(\tau) = \sum_{n \in {\bf Q}} a_{\alpha, n} e^{2\pi \sqrt{-1} n \tau}$$ 
in the obstruction space, the relation
$$\sum_{\alpha \in A_{L}, n<0} a_{\alpha, -n/2}c_{\alpha, n} = 0$$
holds.  In this case the weight $k$ is given by
$$k = \sum_{ \alpha \in A_{L}, n\in {\bf Z}} b_{\alpha, n/2}c_{\alpha, -n}$$
where $b_{\alpha, n}$ are the Fourier coefficients of the Eisenstein series in the obstruction 
space with the constant term
$b_{0, 0} = -1/2$ and $b_{\alpha, 0} = 0$ for $\alpha \not= 0$.}

\medskip

In the following we shall study the divisors $\sum_{\alpha \in A_{L}, n<0} c_{\alpha , n} \calD_{\alpha, n}$
where $c_{\alpha, n}$ depends only on the type of $\alpha$.
Recall that there are 1, 20, 30 and 30 elements in $A_L$ of types
$(00)$, $(0)$, $(1)$ and $(2)$, respectively (see Lemma \ref{types}). 
We denote the vector valued modular form $(f_{\alpha})_{\alpha\in A_L}$ by 
$$(f_{00},\ f_{0},\ f_{1},\ f_{2})$$
where each $f_t$ is the sum of the $f_{\alpha}$ as $\alpha$ varies over the elements of $A_L$ of type $t$.
A calculation shows that the action of a generator $\{ S, T \}$
of $\SL(2,{\bf Z})$ with respect to this basis  is given by
$$\rho^*(T) =
\begin{pmatrix}1&0&0&0
\\0&1&0&0&
\\0&0&\omega^2&0
\\0&0&0&\omega
\end{pmatrix}, \quad 
\rho^*(S) = {-1\over 9}
\begin{pmatrix}1&1&1&1
\\20&-7&2&2
\\30&3&3&-6
\\30&3&-6&3
\end{pmatrix}.
$$

\subsection{Lemma}\label{obs}
{\it The dimension of the space of modular forms of weight $4 = (2+6)/2$ and of 
type $\rho^{*}$ is $2$.  The dimension of
the space of Eisenstein forms of weight $4$ and of type $\rho^{*}$ is also $2$.}

\begin{proof}
The dimension is given by
$$d +dk/12 -\alpha(e^{\pi \sqrt{-1} k/2} \rho^*(S)) -\alpha ((e^{\pi\sqrt{-1}k/3} \rho^*(ST))^{-1}) - \alpha(\rho^*(T))$$
(\cite{B2}, section 4, \cite{F}, Proposition 2.1).
Here $k=4$ is the weight, 
$$d=\dim \{ x \in V : \rho^*(-E)x = (-1)^kx\} =4$$ 
and 
$$\alpha (A)=\sum_{\lambda} \alpha$$ 
where $\lambda$ runs through all eigenvalues of $A$ and $\lambda = e^{2\pi\sqrt{-1}\alpha}$, $0\leq \alpha <1.$
A direct calculation shows that 
$$\alpha(e^{\pi \sqrt{-1} k/2} \rho^*(S)) = 1, \ \alpha ((e^{\pi\sqrt{-1}k/3} \rho^*(ST))^{-1}) = 4/3\ {\rm and}\ 
\alpha(\rho^*(T)) = 1.$$

On the other hand, the space of Eisenstein series is isomorphic to the subspace of $V$ given by
$$\rho^*(T)(x) =x, \ \rho^*(-E)(x) = (-1)^kx$$
(see Remark 2.2 in \cite{F}).
Thus we have the assertion. 
\end{proof}
\smallskip

Next we shall calculate a basis of Eisenstein forms of weight 4 and of type $\rho^{*}$.
Let 
$$E_{1} = G_{4}(\tau, 0,1; 3),\  E_{2} = G_{4}(\tau, 1,0; 3),\ 
E_{3} = G_{4}(\tau,1,1;3),\ E_{4} = G_{4}(\tau,1,2;3)$$ 
be Eisenstein series of weight 4
and level 3 (see \cite{Kob}, Chap. III, \S 3).  Then the action of $S, T$ is as follows: 
$$T: \ E_2 \to E_3 \to E_4 \to E_2,$$
$T$ fixes $E_1$, and $S$ switches $E_1$ and $E_2$, $E_3$ and $E_4$ respectively.

Now we can easily see that a basis of Eisenstein forms of weight 4 and of type $\rho^{*}$ is given by
$$f_{00} = a E_{1} + b(E_{2} + E_{3} + E_{4}),$$
$$f_{0} = (-a-9b) E_{1} + (-3a-7b)(E_{2} + E_{3} + E_{4}),$$
$$f_{1} = (-3a + 3b)(E_{2} + \omega E_{3} + \omega^{2} E_{4}),$$
$$f_{2} = (-3a + 3b)(E_{2} + \omega^{2} E_{3} + \omega E_{4}),$$
where $a, b$ are parameters.
The Fourier expansions of $E_{i}$ are given as follows (see \cite{Kob}, Chap. III, \S 3, Proposition 22):

$$E_{1} =  {(2\pi)^{4} \over 2\cdot 3^{6}}  +  c (-3^{3} q + \cdot \cdot \cdot ),$$

$$E_{2} = c(q^{1/3} + (2^{3} + 1)q^{2/3} + 3^{3} q + \cdot \cdot \cdot ),$$

$$E_{3} = c (\omega  q^{1/3} + (2^{3}+1) \omega^{2} q^{2/3} + 3^{3} q + \cdot \cdot \cdot ),$$

$$E_{4} =  c (\omega^{2}  q^{1/3} + (2^{3}+1) \omega q^{2/3} + 3^{3} q + \cdot \cdot \cdot),$$
where $c = {(-2\pi \sqrt{-1})^{4} \over 3^{4} 3!} = {(2\pi)^4 \over 2\cdot 3^5}$.
Put $a = -{3^{6} \over (2\pi)^{4}}$ and $a = -9b$.  Then

$$f_{00} = -1/2 + 2 \cdot 3^{3} q + \cdot \cdot \cdot ,$$
$$f_{0} = 10 \cdot 3^{3} q + \cdot \cdot \cdot , $$ 
$$f_{1} =  135  q^{2/3} + \cdot \cdot \cdot ,$$
$$f_{2} =  15  q^{1/3} + \cdot \cdot \cdot .$$
It follows from Lemma \ref{obs} that there are no non-zero cusp forms in the obstruction space.
Hence Theorem \ref{freitag} implies that
\subsection{Theorem}\label{}
{\it There exist automorphic forms on $\calD$ of
weight $135$, $15$ with some character whose zero divisors are $(-4/3)$-Heegnear divisor, $(-2/3)$-
Heegnear divisor, respectively}.
\smallskip

Since $(-4/3)$-, $(-2/3)$-Heegnear divisors meet the complex ball with multiplicity 3 (Lemma \ref{restriction}), 
we can take the cube root of these automorphic forms and then we have:

\subsection{Corollary}\label{borcherdspro}
{\it There exist automorphic forms $\Phi_{45}, \Phi_{5}$ on the complex ball $\calB$ of
weight $45$, $5$ whose zero divisors are $(-2)$-Heegnear divisor, $(-1)$-
Heegnear divisor}.

\section{Gritsenko-Borcherds liftings}\label{additive}

In this section, by applying the theory of liftings \cite{B1}, 
we construct a linear system of automorphic forms on $\calB$ with respect to $\Gamma(\sqrt{-3})$
which gives a birational map from the Segre cubic to the Igusa quartic.  

Let $\rho$ be the Weil representation given in \S 4.
A holomorphic map 
$$f : H^+ \to {\bf C}[A_L]$$
 is called a {\it vector valued modular form of weight $k$
and of type} $\rho$ if
$$f(M\tau)=\rho(M)(c\tau+d)^kf(\tau)$$
for any 
$M = \begin{pmatrix}a&b
\\c&d
\end{pmatrix}
\in \SL(2,{\bf Z})$.  
We shall construct a 5-dimensional space of vector valued modular forms of weight 4 
and of type $\rho$.  

Let $V$ be the 5-dimensional subspace of ${\bf C}[A_{L}]$ in Definition \ref{5-dim}.  
First we shall consider the following special vectors $v_{\alpha_0}$ in $V$
(We remark that the following definition of $v_{\alpha_0}$ is similar to the one given in Allcock-Freitag \cite{AF}
to construct liftings in their case).
Let $\alpha_{0}, \alpha_{1}, \alpha_{2}, \alpha_{3}$ be an orthogonal  basis of $A_{L}$ with
$q_L(\alpha_{0}) = -4/3, q_L(\alpha_{1}) = q_L(\alpha_{2}) = q_L(\alpha_{3}) = -2/3$.
If we take $\alpha_{0}$, then such basis is uniquely determined up to signs.
For each such basis, we define a vector 
 $v_{\alpha_0} = (c_{\alpha})_{\alpha \in A_{L}}$
in ${\bf C}[A_{L}]$ as follows:
$$c_{\alpha} = 1, 0, -1$$
according to 
$$\prod_{i} \langle \alpha, \alpha_{i}  \rangle  = 1, 0, -1 \in {\bf F}_{3}.$$
For example, 
assume $\alpha_{0} = (1,0,0,0)$,  $\alpha_{1} = (0,1,0,0)$, $\alpha_{2} = (0,0,1,0)$,
$\alpha_{3} = (0,0,0,1) \in A_L=({\bf F}_3)^4$.
Then $c_{\alpha} \not= 0$ if and only if $\alpha \in \{ (\pm 1, \pm 1, \pm 1, \pm 1) \}$.

\subsection{Lemma}\label{invariant}
{\it Let $v_{\alpha_0} = (c_{\alpha})$ be as above.  Then $\rho (S)(v_{\alpha_0}) = v_{\alpha_0}$, 
$\rho (T)(v_{\alpha_0}) = \omega v_{\alpha_0}$.   
Moreover $r_{\alpha_{i}}(v_{\alpha_0}) = -v_{\alpha_0}$ for the reflection $r_{\alpha_{i}}$
associated with $\alpha_{i}$.  } 
\begin{proof}
It suffices to prove the case $\alpha_{0} = (1,0,0,0)$,  $\alpha_{1} = (0,1,0,0)$, $\alpha_{2} = (0,0,1,0)$,
$\alpha_{3} = (0,0,0,1)$.  Let $M = \{ (\pm 1, \pm 1, \pm 1, \pm 1) \}$.
Then $v_{\alpha_0} = \sum_{\alpha \in M} c_{\alpha} e_{\alpha}.$  If $c_{\alpha} \not= 0$, then $q_L(\alpha) = -4/3$.  
Hence $\rho (T)(v_{\alpha_0}) = \omega v_{\alpha_0}$.
Next consider 
$$\rho (S) v_{\alpha_0} = -{1\over 9} \sum_{\beta \in A_L} (\sum_{\alpha\in M} c_{\alpha} e^{-2\pi \sqrt{-1} \la \alpha, \beta \ra}) e_{\beta}.$$
A direct calculation shows that the coefficient 
$$\sum_{\alpha\in M} c_{\alpha} e^{-2\pi \sqrt{-1} \la \alpha, \beta \ra}$$
of $e_{\beta}$ is $0$ if $\beta \notin M$, $9$ if $\beta \in M$, $c_{\beta} = -1$, and $-9$ if $\beta \in M$, $c_{\beta} =1$.
The last assertion follows from the definition of $v_{\alpha_0}$.
\end{proof}

It follows from Lemma \ref{invariant} that $v_{\alpha_0}$ is contained in $V$.  Thus we have fifteen elements $v_{\alpha_0}$
in $V$ where $\alpha_0$ is fifteen $(-4/3)$-vectors in $A_L/\{ \pm 1\}$.

\subsection{Lemma}\label{} 
{\it As a $\O(q_L)$ module, $V$ is irreducible.}
\begin{proof}
If $W$ is a irreducible representation of ${\mathfrak S}_6$ and $\dim \ W \geq 2$, then $\dim \ W \geq 5$.  Hence it suffices to see that
there are no 1-dimensional subspaces invariant under the action of ${\mathfrak S}_6$.  If such 1-dimensional subspace exists, 
then all vectors in $V$ are invariant under the action of ${\mathfrak S}_6$.
However any special vector $v_{\alpha_0}$ as above is not invariant under the action of ${\mathfrak S}_6$.  This is a contradiction.
\end{proof}

\smallskip

Let $\eta(\tau)$ be the Dedekind eta function.  Then
$$\eta (\tau + 1)^{8}= \omega \cdot \eta(\tau)^{8},$$
$$\eta(-1/\tau)^{8} = \tau^{4} \cdot \eta(\tau).$$
Therefore, for $v\in V$, $\eta(\tau)^8 \cdot v = (\eta(\tau)^8 \cdot c_{\alpha})_{\alpha \in A_{L}}$ is a vector valued modular form of weight
4 and of type $\rho$.  By applying Gritsenko-Borcherds lifting (\cite{B1}, Theorem 14.3), we have

\subsection{Lemma}\label{5dim}
{\it There is a $5$-dimensional space of holomorphic automorphic forms of weight $6$
on $\calD$ with respect to $\tilde{\O}(L)$ on which $\O(q_{L})$ acts irreducibly.}
\begin{proof}
It suffices to see that the lifting of $\eta(\tau)^8 v$ is non-zero.  Then the assertion follows from the Schur's lemma. 
We use Theorem 14.3 in \cite{B1}.  
We consider the orthogonal decomposition of
$L = U\oplus U(3) \oplus A_2\oplus A_2$.  
Let $z, z'$ be a basis of $U$ with $z^2 =z'^2 =0, \la z, z'\ra = 1$, and let $K = z^{\perp}/{\bf Z} z = U(3)\oplus A_2\oplus A_2
\subset L$.  Let $e, f$ be a basis of $U(3)$ with $e^2=f^2=0, \la e, f\ra =3$.  We consider the Fourier expansion around $z$.
Since $\eta(\tau)^8 = q^{1/3} + \cdot \cdot \cdot$, the initial term of 
$$\eta(\tau)^8 v = \sum_{\alpha \in A_L} e_{\alpha} \sum_{n \in {\bf Q}} c_{\alpha}(n) e^{2\pi\sqrt{-1}n\tau}$$  
is 
$$\sum_{\alpha \in A_L, \alpha^2 =2/3} c_{\alpha} q^{1/3}.$$  If we take $\lambda = (e+f)/3$, then 
$\la \lambda, \lambda \ra = 2/3 > 0$ and hence $\lambda$ has positive inner products with all elements in the interior of
the Weyl chamber.  Also note that $L^*/L = K^*/K$.  We choose $v = (c_{\alpha}) \in V$ satisfying $c_{\lambda} \not=0$.
Now it follows from \cite{B1}, Theorem 14.3 that the Fourier coefficient of $e^{2\pi \sqrt{-1} \la \lambda, Z\ra }$ in 
the lifting of $\eta(\tau)^8 v$ is equal to 
$$c_{\lambda}(\lambda^2/2)\cdot e^{2\pi \sqrt{-1} \la \lambda, z'\ra}= c_{\lambda}(1/3) = c_{\lambda}.$$
Hence the lifting of $\eta(\tau)^8 v$ is non-zero.
\end{proof}

Let $\alpha_0$ be a $(-4/3)$-vector in $A_L$ and let $v_{\alpha_0}$ be the element in $V$ as above.
Let $F_{\alpha_0}$ be the restriction of the Gritsenko-Borcherds lifting of $\eta(\tau)^8 \cdot v_{\alpha_0}$ to the complex
ball $\calB$.
Then 

\subsection{Theorem}\label{}
{\it $F_{\alpha_0}$ is a holomorphic automorphic form of weight $6$
on $\calB$ with respect to $\Gamma(\sqrt{-3})$ which vanishes exactly along the $(-2)$-Heegner divisor ${\calH}_{\alpha_0}$ with multiplicity one
and the $(-1)$-Heegner divisors ${\calH}_{\alpha_1}, {\calH}_{\alpha_2}, {\calH}_{\alpha_3}$ with multiplicity three.}
\begin{proof}
First recall that the reflection $r_{\alpha}$ is induced from the reflection $r_{a, -1}$ of $\Lambda$ where 
$a \in \Lambda$ is a short or long root with $a\ \mod \ \sqrt{-3}\Lambda = \alpha$.
It follows from Lemma \ref{invariant} and the $\O(q_L)$-equivariance of the lifting that $F_{\alpha_0}$ vanishes along 
${\calH}_{\alpha_i}$ $(i=0,1,2,3)$.  Moreover the trireflection $r_{a,\omega}$ associated to a short root $a$ is contained in 
$\Gamma(\sqrt{-3})$, $F_{\alpha_0}$ vanishes along ${\calH}_{\alpha_i}$ $(i=1,2,3)$ with multiplicity 3.
Then the product of 15 $F_{\alpha_0}$ has weight $90$ and vanishes along $(-2)$-Heegner divisors with at least
multiplicity one and along $(-1)$-Heegner divisors with at least multiplicity ${3\cdot 15 \cdot 3 \over 15} =9$.
On the other hand $\Phi_{45} \cdot \Phi_{5}^9$ has weight $90$ and vanishes 
along $(-2)$-Heegnear divisors with exactly multiplicity one and along $(-1)$-Heegnear divisors with multiplicity $9$
(Corollary \ref{borcherdspro}).
Then the ratio $\prod_v F_v / \Phi_{45} \cdot \Phi_{5}^9$ has weight zero and holomorphic, and hence it is constant
by Koecher principle.
\end{proof}

Since trireflections are contained in $\Gamma (\sqrt{-3})$, the covering ${\calB} \to {\calB}/\Gamma(\sqrt{-3})$ is ramified along $(-2/3)$-Heegnear divisors.  Hence we have

\subsection{Theorem}\label{}
{\it The zero divisor $(F_{\alpha_0})$ on $\bar{\calB}/\Gamma(\sqrt{-3})$ is
$\bar{\calH}_{\alpha_0} + \bar{\calH}_{\alpha_1} + \bar{\calH}_{\alpha_2} + \bar{\calH}_{\alpha_3}$.}

\smallskip

\subsection{Lemma}\label{Cayley} 
{\it Let $\alpha \in A_L$ with $q_L(\alpha)=-4/3$.  Then the Heegner divisor $\bar{\calH}_{\alpha}$ coincides with a Cayley cubic on $X$.}

\begin{proof}
The zero divisor of $\Phi_{45}$ on $\bar{\calB}/\Gamma(\sqrt{-3})$ is the union of $\bar{\calH}_{\alpha}$ where $\alpha$ varies over
15 $(-4/3)$-vectors in $A_L/\{ \pm 1\}$.  On the other hand, as mentioned as above, 
the covering ${\calB} \to {\calB}/\Gamma(\sqrt{-3})$ is ramified along $(-1)$-Heegnear divisors.  Hence
the zero divisors of $\Phi_5^3$ on $\bar{\calB}/\Gamma(\sqrt{-3})$ is the union of  $\bar{\calH}_{\alpha}$ where
$\alpha$ varies over 15 $(-2/3)$-vectors in $A_L/\{ \pm 1\}$.  
If $\alpha_1, \alpha_2, \alpha_3$ are mutually orthogonal $(-2/3)$-vectors, then
$\bar{\calH}_{\alpha_1} + \bar{\calH}_{\alpha_2} + \bar{\calH}_{\alpha_3}$ is the union of three planes on $X$ (Proposition \ref{15planes}), that is, 
a linear section of $X$ in ${\bf P}^4$.  By comparing the weights of $\Phi_{45}$ and $\Phi_{5}^3$, we can see that 
each $\bar{\calH}_{\alpha}$ with $q_L(\alpha) = -4/3$ is also a linear section, that is, a cubic surface.
Since $\bar{\calH}_\alpha$ with $q_L(\alpha) = -4/3$ contains four cusps, it should be isomorphic to a Cayley cubic on $X$.
\end{proof}

\noindent
Hence we conclude:

\subsection{Lemma}\label{} 
{\it The divisor $(F_{\alpha_0}) = \bar{\calH}_{\alpha_0} + \bar{\calH}_{\alpha_1} + \bar{\calH}_{\alpha_2} + \bar{\calH}_{\alpha_3}$ is
a quadric section of the Segre cubic $X \subset {\bf P}^4$ where
 $\bar{\calH}_{\alpha_0}$ is a Cayley cubic and $\bar{\calH}_{\alpha_i} \ (i=1,2,3)$ are planes.}

\smallskip
For an orthogonal basis $\{ \alpha_0, \alpha_1, \alpha_2, \alpha_3\}$ of $A_L$, we can easily see that any isotropic vector in $A_L$ is
perpendicular to $\alpha_i$ for some $i$ (see Lemma \ref{types}).  Hence the divisor 
$\bar{\calH}_{\alpha_0} + \bar{\calH}_{\alpha_1} + \bar{\calH}_{\alpha_2} + \bar{\calH}_{\alpha_3}$ contains 10
nodes of $X$.  Thus
each $(F_{\alpha_0})$ passes the ten nodes of $X$, and hence
the five dimensional linear system of automorphic forms has the ten nodes as base points.  The linear system defines a rational map
$\varphi : X - - \to {\bf P}^5$.
\smallskip

\subsection{Theorem}\label{}
{\it The image of $\varphi$ is the Igusa quartic, that is the dual of $X$.}
\begin{proof}
This follows from \cite{H2}, Theorem 3.3.12.
\end{proof}


\end{document}